\documentclass[12pt]{article}
\input epsf      
\usepackage{graphicx}% Include figure files
\usepackage{amssymb}
\usepackage{amsfonts}
\usepackage{amsthm}
\usepackage{amsfonts}

\newcommand{\reals}{\mathbb{R}}

\newcommand{\half}{\mathbb{H}}

\newcommand{\compose}{\circ}

\newcommand{\ba}[1]{\begin{array}{#1}}
\newcommand{\ea}{\end{array}}

\newcommand{\be}{\begin{equation}}
\newcommand{\ee}{\end{equation}}
\newcommand{\bea}{\begin{eqnarray}}
\newcommand{\eea}{\end{eqnarray}}
\newcommand{\beann}{\begin{eqnarray*}}
\newcommand{\eeann}{\end{eqnarray*}}

\def\reff#1{(\ref{#1})}
\begin{document}

\bibstyle{ams}

\title{Monte Carlo comparisons of the self-avoiding walk and SLE
as parameterized curves } 

\author{Tom Kennedy
\\Department of Mathematics
\\University of Arizona
\\Tucson, AZ 85721
\\ email: tgk@math.arizona.edu
}
\maketitle

\begin{abstract}

The scaling limit of the two-dimensional self-avoiding walk (SAW) is believed
to be given by the Schramm-Loewner evolution (SLE) with the parameter $\kappa$
equal to $8/3$. The scaling limit of the SAW has a natural parameterization 
and SLE has a standard parameterization using the half-plane capacity. 
These two parameterizations do not correspond with one another. To make 
the scaling limit of the SAW and SLE agree as parameterized curves,
we must reparameterize one of them. 
We present Monte Carlo results that show that if we reparameterize
the SAW using the half-plane capacity, then it agrees well with SLE with its
standard parameterization.
We then consider how to reparameterize SLE to make 
it agree with the SAW with its natural parameterization. We argue using
Monte Carlo results that the so-called $p$-variation of the SLE curve with
$p=1/\nu=4/3$ provides a parameterization that corresponds to the natural
parameterization of the SAW. 

\end{abstract}

\bigskip

\newpage

\section{Introduction}
\label{intro}

The scaling limit of the two dimensional self-avoiding walk (SAW)
is believed to be the Schramm-Loewner evolution (SLE) 
with $\kappa=8/3$ \cite{lsw_saw}. 
Previous Monte Carlo studies have compared the scaling limit of the SAW and
SLE$_{8/3}$ as unparameterized curves \cite{tk_saw_sle_one,tk_saw_sle_two}. 
These studies considered random variables that did not depend on how
the curve were parameterized and found excellent agreement between
the distributions for the SAW and SLE$_{8/3}$. 
The goal of this paper is to use Monte Carlo simulations 
to understand how the SAW and SLE$_{8/3}$ curves 
should be parameterized so that they agree as parameterized curves. 

We give a brief definition of the SAW. For a detailed treatment we refer 
the reader to \cite{ms}. The relationship of the SAW to SLE$_{8/3}$ 
is the subject of \cite{lsw_saw}.
We restrict our attention to the SAW on a two-dimensional lattice 
restricted to the upper half plane. Fix a positive integer $N$. 
We consider all nearest neighbor walks which begin at the origin, 
have exactly $N$ steps, remain in the upper half plane (except for
the starting point) and never visit the same site more than once.  
There are a finite number of such walks, so we can put the uniform probability
measure on this set of walks. We then attempt to construct the scaling limit
by the following double limit. We first let $N \rightarrow \infty$.
For the case of the half-plane this limit was proved to exist in 
\cite{lsw_saw}, but there is no general proof of its existence.
We then take the lattice spacing to zero. This limit has not been proved 
to exist, but it is believed that it does and gives a measure on 
curves in the upper half plane which start at the origin and go to $\infty$. 
The SAW curves are simple (do not intersect themselves) by definition.
This does not insure that in the scaling limit the measure is
supported on simple curves, but it is believed that it is. 

SLE was introduced in \cite{schramm} and studied further in \cite{rs}.
Expositions may be found in \cite{lawler} and \cite{werner}.
In the case of $\kappa=8/3$, SLE satisfies the restriction property and 
so belongs to the one parameter family of restriction measures 
studied by Lawler, Schramm and Werner \cite{lsw_restrict}. 
They showed that the probabilities of certain 
events for SLE$_{8/3}$ may be computed from certain conformal maps. In 
some cases they may be computed explicitly. Past Monte Carlo simulations
have tested the equivalence of the scaling limit of the SAW and SLE$_{8/3}$
by computing these probabilities for the SAW by simulation and comparing 
with the exact result for SLE$_{8/3}$. Here are two examples.
Let $P$ be a point on the horizontal axis and consider the minimum distance 
of the SLE$_{8/3}$ curve to the point. 
This is a random variable whose distribution
may be computed explicitly for SLE$_{8/3}$. For the second example 
consider a vertical line given by $x=c$ for some $c \neq 0$. 
The SLE$_{8/3}$ curve intersects it. Let $Y$ be the distance from 
the horizontal axis to the lowest intersection. The distribution of $Y$ 
may be explicitly computed for SLE$_{8/3}$. 
Note that for both of these random variables the time parameterization
of the curve plays no role. For these two examples and 
several others, excellent agreement was found between the SLE$_{8/3}$
exact results and the simulation results for the scaling limit of the SAW
\cite{tk_saw_sle_one,tk_saw_sle_two}. 

The SLE curve is usually parameterized so that its half-plane capacity 
grows linearly with the time parameter $t$. Before we take the scaling
limit, the SAW has a natural parameterization given by the number of steps, 
or equivalently the distance along the curve since all steps are 
nearest neighbor. With a suitable scaling this should give a natural 
parameterization of the scaling limit of the SAW. 
Equipped with these natural parameterizations, the scaling limit of 
the SAW and the SLE$_{8/3}$ parameterized curves are not the same. 
We illustrate this with a simulation in section \ref{need}. 
How should these curves be reparameterized so that they agree?

There are two ways to ask this question. We can ask how to 
reparameterize the scaling limit of the SAW to make it agree 
with SLE$_{8/3}$ with its parameterization using half-plane capacity.
The answer is obvious - 
one should use the half-plane capacity of the scaling limit of
the SAW to parameterize it. We study this in section \ref{reparamSAW}.

A less obvious question is to ask how we should reparameterize 
SLE$_{8/3}$ to make it agree with the scaling limit of the SAW  
with its natural parameterization.
In section \ref{fvar} we give simulation results that 
indicate that for the SAW with its natural parameterization, 
the $p$-variation of the walk exists and is a deterministic, linear
function of time if we take $p=1/\nu=4/3$. In other words, one can 
recover the natural parameterization of the SAW by computing its 
$1/\nu$ variation. 
Some care is needed here since different definitions of the 
variation give different values. However, it appears they differ
only by a multiplicative constant. 
In section \ref{reparamSLE} we use the
$1/\nu$ variation of the SLE to reparameterize the SLE curve and 
compare this with the SAW with its natural parameterization.
Section \ref{simulating} discusses various computational aspects 
of our simulations of the SAW and the SLE. 

\section{Need for a random reparameterization}
\label{need}

We review some facts about SLE. We refer the reader to \cite{lawler} for 
more detail. We restrict our attention to chordal SLE in the upper half plane.
It is defined by the differential equation
\be
\dot{g}_t(z) = {2 \over g_t(z) - \sqrt{\kappa} B_t}
\label{sle_eq}
\ee
and the initial condition $g_0(z)=z$. 
The dot denotes differentiation with respect to time $t$, 
$B_t$ is a standard real-valued Brownian motion, and
$z$ belongs to the upper half plane $\half$.  
For some initial points $z$, the solution only exists for a finite
time interval. The set of initial points for which the solution 
no longer exists at time $t$ is a random subset $K_t$ of $\half$. 
For $\kappa \le 4$, $K_t$ is just a simple curve $\gamma(t)$, i.e., 
a curve that does not intersect itself.  
We use $\gamma[0,t]$ to denote the image of the curve for the 
time interval $[0,t]$. So $\gamma[0,t]=K_t$. 

The parameterization of the curve $\gamma$ that results from 
eq. \reff{sle_eq} corresponds to the half-plane capacity of the 
curve which is defined as follows. 
Since $\gamma$ is simple, $\half \setminus \gamma[0,t]$
is simply connected. So there is a conformal map of this domain onto $\half$. 
If we require the map to send $\infty$ to itself, then 
this map has an expansion about $\infty$ of the form
\be
g(z)= c z + a_0 + \sum_{n=1}^\infty a_n z^{-n}
\ee
The map is uniquely determined if we require $c=1$ and $a_0=0$. 
The half plane capacity of $\gamma[0,t]$ is then the coefficient $a_1$. 
We denote it by $hcap(\gamma[0,t])$.

The above definition applies to any simple curve $\gamma$. For the 
SLE curve we have 
\be
hcap(\gamma[0,t])  = 2t
\ee
If we were to parameterize SLE so that 
\be
hcap(\gamma[0,t])=b(t)
\ee
for some differentiable, increasing function $b(t)$, then the $2$ in 
the Schramm-Loewner equation would be replaced by $\dot{b}(t)$.

The scaling property of Brownian motion leads to a scaling property for 
SLE: $\gamma[0,t]$ has the same distribution as $\sqrt{t} \, \gamma[0,1]$. 
(Note that if we rescale a set by $r$, its half-plane capacity scales 
by a factor of $r^2$.) 
Thus 
\be
E[ \gamma(t)^2] = c \, t
\ee
where $c=E[\gamma(1)^2]$.

The SAW has a natural parameterization. 
Let $W(n)$ be the infinite SAW in the upper half plane on a lattice 
with unit lattice spacing. We extend $W(t)$ to all $t \ge 0$ by 
linearly interpolating between the integer times.
The average distance from the origin to $W(n)$ 
is believed to grow like $n^\nu$ with $\nu=3/4$.
The limit of taking the lattice 
spacing to zero can then be obtained by defining 
\be
\omega(t) = \lim_{n \rightarrow \infty} n^{-\nu} W(nt)
\label{scalinglimit}
\ee
This limit is believed to exist and give a probability measure on 
simple curves in the upper half plane starting at the origin, but this
has not been proved. 
This definition of the scaling limit is analogous to Brownian motion
as the scaling limit of an ordinary random walk. If we take $W(n)$ to be
an ordinary random walk in the plane and $\nu=1/2$, then \reff{scalinglimit}
would converge to two-dimensional Brownian motion.
We refer to this parameterization of the scaling limit of 
the SAW as its natural parameterization.
If this limit exists, then $\omega(t)$ has the same distribution 
as $t^\nu \omega(1)$. Thus 
\be
E [\omega(t)^2] = b \, t^{2 \nu}
\ee
where $b=E [\omega(1)^2]$.

The behavior of $E[\gamma(t)^2]$ and $E[\omega(t)^2]$ show that 
to make the scaling limit of the SAW and SLE$_{8/3}$ 
agree we must at least do a deterministic reparameterization. 
If there is a reparameterization, $\hat{\gamma}(t)=\gamma(\phi(t))$,
of SLE$_{8/3}$ with a deterministic $\phi(t)$ that makes $\hat{\gamma}(t)$ 
and $\omega(t)$ agree in distribution, then we must have 
$ E[\hat{\gamma}(t)^2] = E[\omega(t)^2]$ and so the reparameterization must 
be 
\be
\hat{\gamma}(t) = \gamma(a \, t^{2 \nu})
\ee
with $a=b/c$. 

\begin{figure}[tbh]
\includegraphics{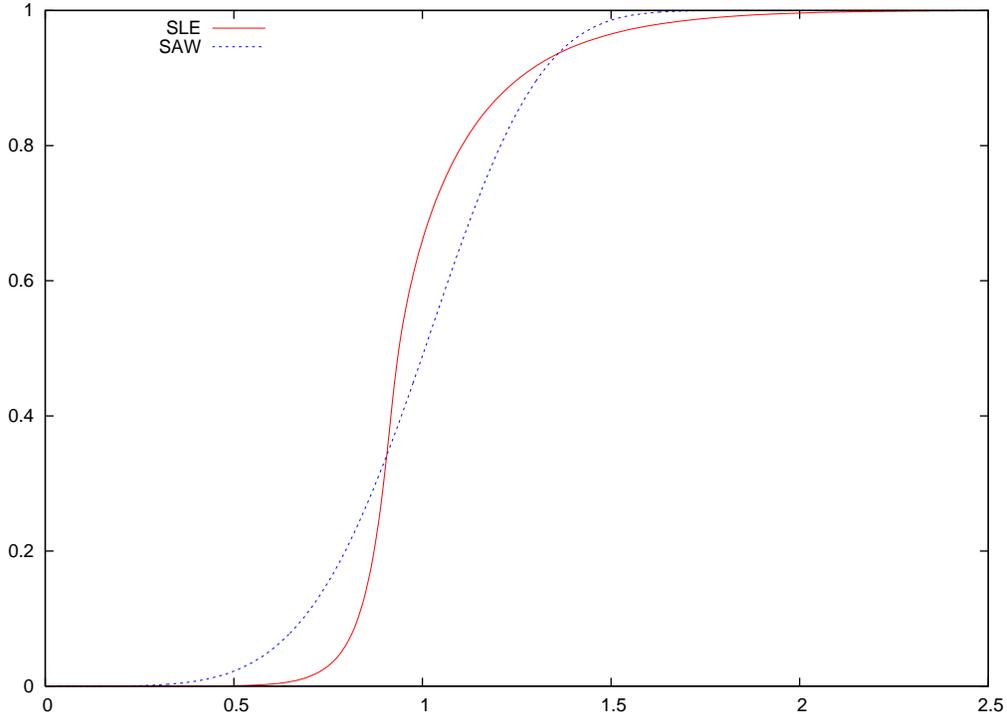}
\caption{ The distribution, $P(R \le t)$, for the distance $R$ from the 
origin to $\gamma(1)$ (labeled by SLE) and to $\omega(1)$ (labeled by 
SAW). Both random variables have been rescaled so that they have mean one. 
}
\label{saw_sle_fix_time_dist}
\end{figure}

\begin{figure}[tbh]
\includegraphics{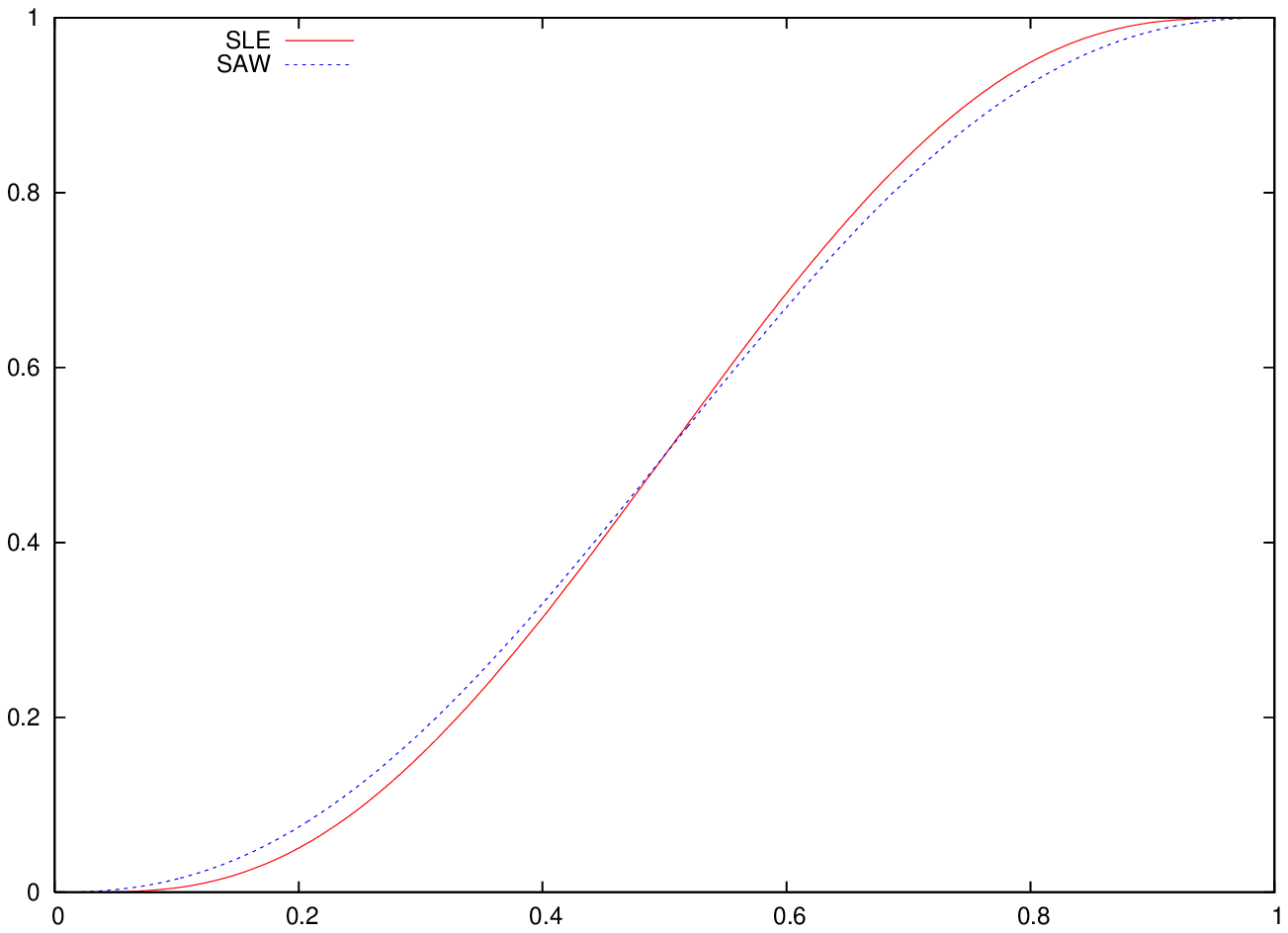}
\caption{ The distribution, $P(\Theta \le t)$, for the polar angle $\Theta$ 
of the points $\hat{\gamma}(1)$ (labeled by SLE) and $\omega(1)$ (labeled by 
SAW). 
}
\label{saw_sle_fix_time_angle}
\end{figure}

It is easy to see from simulations that this deterministic reparameterization 
does not work, i.e., $\hat{\gamma}(t)$ and $\omega(t)$ do not have the 
same distribution as parameterized curves. To show this we consider 
$\hat{\gamma}(1)$ and $\omega(1)$. If they have the same distribution, then 
$\gamma(1)/E[|\gamma(1)|]$ and $\omega(1)/E[|\omega(1)|]$ have 
the same distribution. 
For each of these two random points, we 
compute two random variables - the distance $R$ from the origin to the random 
point and the polar angle $\Theta$ of the random point.
Our rescaling means that $E[R]=1$ for both the SLE and the SAW. 

The distributions of $R$ for the SLE and the SAW are shown 
in figure \ref{saw_sle_fix_time_dist}. The distributions of $\Theta$ are
shown in figure \ref{saw_sle_fix_time_angle}. 
We do not show any error bars in the figures. The width of the error bars
would be hard to see on this scale.
Clearly the distributions are different.

Throughout this paper our plots of the 
distributions of random variables are plots of their cumulative
distribution functions rather than their densities. The cumulative
distribution is the function that is actually computed in the simulation.
Computing the density would require taking a numerical derivative.

\section{Reparameterizing the SAW}
\label{reparamSAW}

In this section we reparameterize the scaling limit of the SAW 
so that it agrees with SLE with its usual parameterization. 
For $t>0$ we define $C_t$ to be the random time for the scaling limit
of the SAW when
\be
hcap(\omega[0,C_t])=2t
\ee
We define $\hat{\omega}(t)=\omega(C_t)$ so that 
$hcap(\hat{\omega}[0,t])=2t$. Then $\hat{\omega}(t)$ and $\gamma(t)$ 
should have the same distribution as parameterized curves. 
We test this by comparing the distributions of $R$ and $\Theta$ for 
$\gamma(1)$ and $\hat{\omega}(1)=\omega(C_1)$. As before, $R$ is the distance
from the random point to the origin and $\Theta$ is the polar angle of 
the random point. 

\begin{figure}[tbh]
\includegraphics{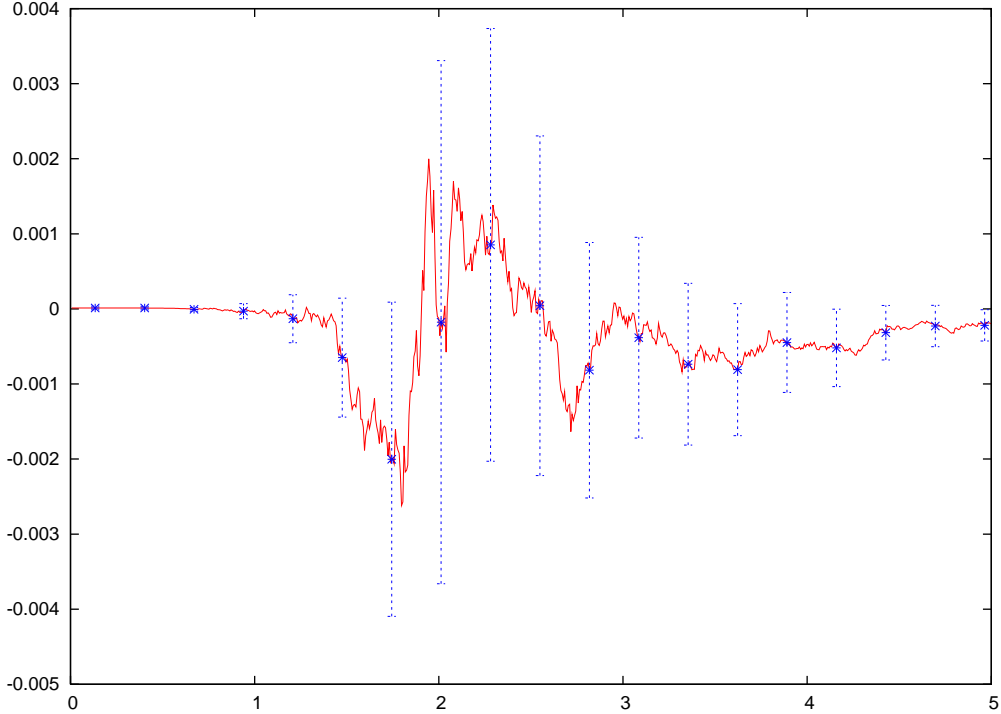}
\caption{The difference of the distributions of $R$ for $\hat{\omega}(1)$,
the SAW parameterized using capacity, and $\gamma(1)$, the SLE with its
usual parameterization.}
\label{saw_sle_fix_cap_dist}
\end{figure}

If we plot the distributions of $R$ for $\gamma(1)$ and 
$\hat{\omega}(1)$, the two curves are virtually indistinguishable. 
So in figure \ref{saw_sle_fix_cap_dist} we plot the difference of 
these two distributions.
The important thing to note about this figure is the scale on the 
vertical axis. The difference between the two distributions is typically
on the order of a tenth of a percent.
We caution the reader that in figure \ref{saw_sle_fix_time_dist}
the random variable $R$ was rescaled so that its mean was 1. 
This is not done in figure \ref{saw_sle_fix_cap_dist}. The mean of $R$ 
here is slightly greater than $2$. In particular, the sharp spikes in 
figure \ref{saw_sle_fix_cap_dist} just left of $R=2$ correspond to where 
the distribution is increasing sharply. In figure \ref{saw_sle_fix_time_dist}
this sharp increase occurs just left of $R=1$. 

\begin{figure}[tbh]
\includegraphics{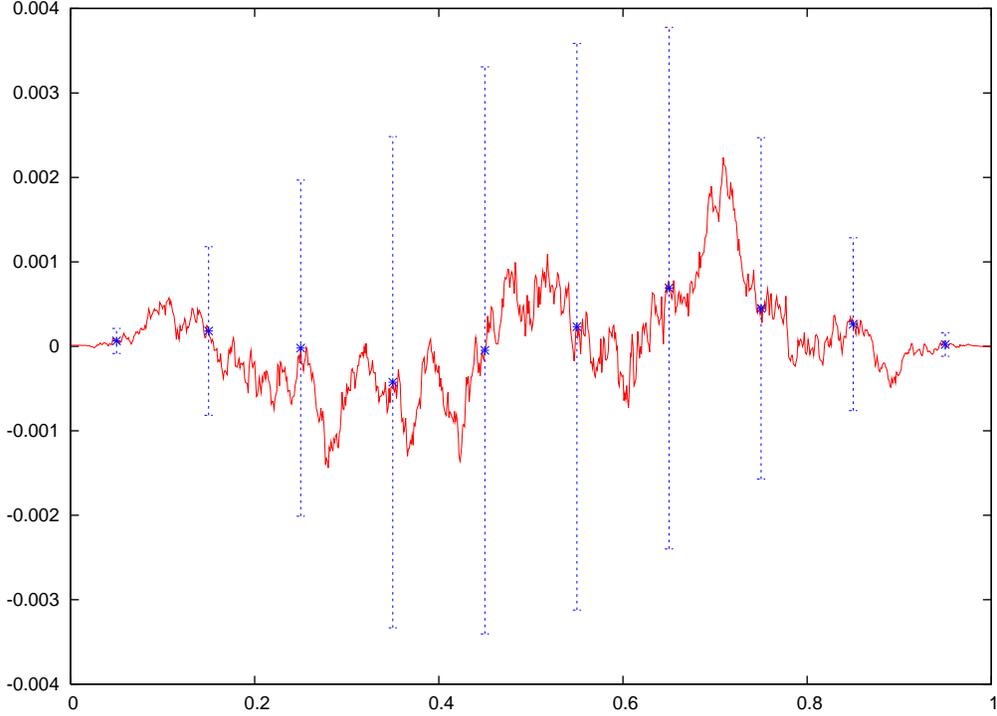}
\caption{The difference of the distributions of $\Theta$ for $\hat{\omega}(1)$,
the SAW parameterized using capacity, and $\gamma(1)$, the SLE with its
usual parameterization.}
\label{saw_sle_fix_cap_angle}
\end{figure}

The difference of the distributions of $\Theta$ 
for $\gamma(1)$ and $\hat{\omega}(1)$ is shown in 
\ref{saw_sle_fix_cap_angle}. Again, the most important feature 
of this plot is the scale on the vertical axis. These two distributions 
also differ by on the order of a tenth of a percent.

The error bars in figures \ref{saw_sle_fix_cap_dist} and 
\ref{saw_sle_fix_cap_angle} and subsequent figures are only statistical
errors, i.e., the error arising from not running the Monte Carlo 
simulations forever. The error bars shown are two standard deviations.
The standard deviation is estimated using the techniques of batched
means. The simulations of the SAW and the SLE both 
contain systematic errors. These are discussed in section 
\ref{simulating}.

\bigskip

To further test the equivalence of $\gamma(t)$ and $\hat{\omega}(t)$, 
we consider the covariances of these two processes. We think of 
these processes as taking values in $\reals^2$ rather than the complex
plane and let $\gamma_i(t)$ and $\hat{\omega}_i(t)$ denote their coordinates 
with $i=1,2$. The covariance of SLE is given by 
\be
C_{ij}(s,t) = E \, [\gamma_i(s) \gamma_j(t)] 
- E \, [\gamma_i(s)]  \, E \, [\gamma_j(t)] 
\ee
By scaling, the dependence of the last term on $s$ and $t$ is trivial:
\be
E \, [\gamma_i(s)]  \, E \, [\gamma_j(t)] = \sqrt{st} \, 
E \, [\gamma_i(1)]  \, E \, [\gamma_j(1)]
\ee
So we will only study $E \, [\gamma_i(s) \gamma_j(t)]$.
By scaling, 
\be
E \, [\gamma_i(s) \gamma_j(t)] = s \, E \, [\gamma_i(1) \gamma_j(t/s)] 
= t \, E \, [\gamma_i(s/t) \gamma_j(1)] 
\ee
So it suffices to study this quantity with $s=1$ and $0 \le t \le 1$. 
We can rewrite this as
\be
 E  \, [\gamma_i(1) \gamma_j(t)] 
=E  \, [\gamma_i(t) \gamma_j(t)] 
+ E  \, [(\gamma_i(1) - \gamma_i(t)) \gamma_j(t)] 
\ee
By scaling, the first term is just $c_{ij} \, t$ where 
$c_{ij}= E  \, [\gamma_i(1) \gamma_j(1)]$.  
So we simulate just the second term. So let
\be
\rho_{ij}(t)= E  \, [(\gamma_i(1) - \gamma_i(t)) \gamma_j(t)] 
\ee
For the SAW parameterized by capacity, we define 
\be
\hat{\sigma}_{ij}(t)= E  \, [(\hat{\omega}_i(1) - \hat{\omega}_i(t)) 
\hat{\omega}_j(t)] 
\ee

\begin{figure}[tbh]
\includegraphics{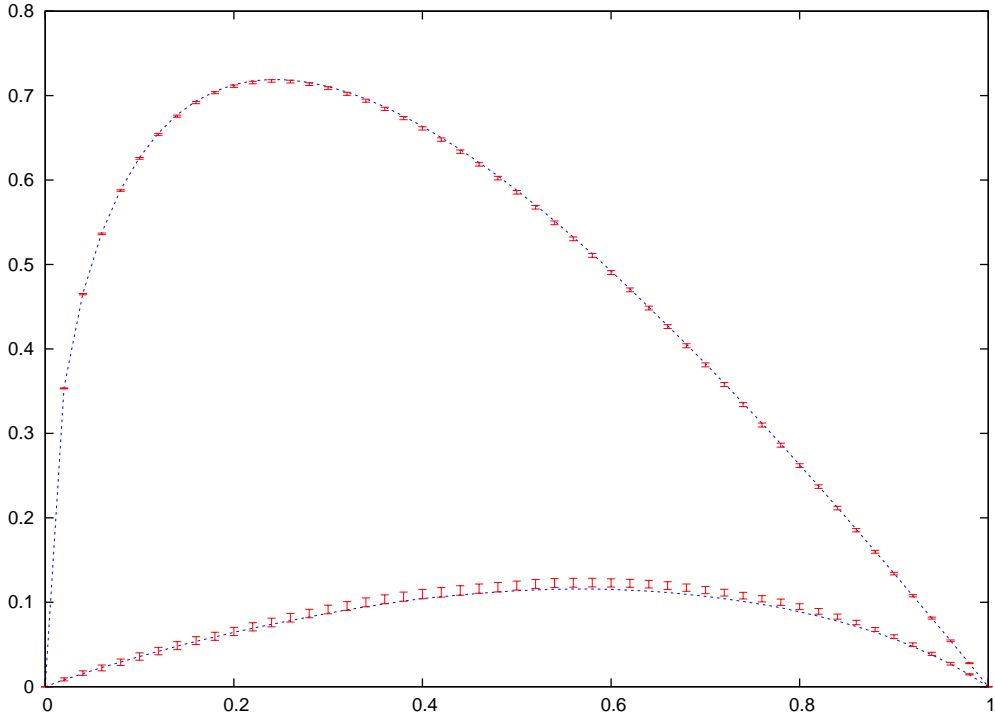}
\caption{ The covariance $\rho_{ij}(t)$ for SLE and 
the covariance $\hat{\sigma}_{ij}(t)$ for the SAW using capacity as its
parameterization. 
The top curve is $i=j=2$, the bottom curve is $i=j=1$. 
The SLE covariance is plotted with a curve and no error bars, while 
individual points with error bars are plotted for the SAW. 
}
\label{saw_sle_fix_cap_covar}
\end{figure}

The functions $\rho_{ij}(t)$ and $\hat{\sigma}_{ij}(t)$ are shown 
in figure \ref{saw_sle_fix_cap_covar} for $i=j=1$ and $i=j=2$. 
(It is not hard to see that by symmetry these functions are zero 
when $i \ne j$.)
The functions agree so well that we have plotted individual points 
for $\hat{\sigma}_{ij}(t)$ and an interpolating curve for $\rho_{ij}(t)$ 
so that the two functions may be distinguished.
It is much harder to simulate the SAW at a constant capacity than to 
simulate the SLE. Consequently the statistical errors for the SAW are 
larger. So we have only shown error bars for the SAW points. 
The difference between the covariances of SLE and the SAW parameterized 
by capacity is small, but appears to be greater than the statistical
error shown by the error bars in the figure. We emphasize again 
that these error bars do not include the systematic errors. 
Computing the capacity of a SAW is time consuming. This limits 
the SAW simulation to walks of modest length - one of the sources of 
systematic error.

\section{$p$-variation of the SAW}
\label{fvar}

We now consider the question of how to reparameterize SLE$_{8/3}$ so that it
agrees with the scaling limit of the SAW with its natural parameterization.
If one plots points on a SAW that are equally spaced in time, then the 
points appear to be equally spaced along the SAW path. In other words, 
the natural parameterization of the SAW corresponds to using the length 
of the walk as the parameter. Of course, this is only a heuristic 
statement since the SAW is not expected to have finite variation. 
However, we will see that the so called ``$p$-variation'' provides 
a way to make sense of the length along the SAW. 

If $X(s)$ is a stochastic process, then the $p$-variation is defined 
as follows. Let $0=t_0^n < t_1^n < t_2^n \cdots < t_{k_n}^n = t$ be a
sequence of partitions of $[0,t]$. Let $\Pi_n$ denote the $n$th partition
and let $||\Pi_n||$ be the width of the largest subinterval of $\Pi_n$. 
We assume that $||\Pi_n||$ goes to zero as $n$ goes to infinity. 
We define 
\be 
var_p((X(s))_{0 \le s \le t},\Pi_n)=
\sum_{j=1}^{k_n} |X(t^n_j)-X(t^n_{j-1})|^p
\label{fvar_def}
\ee
With some abuse of the notation we will write 
$var_p((X(s))_{0 \le s \le t},\Pi_n)$ as just $var_p(X(0,t),\Pi_n)$.
We define the $p$-variation to be the limit as the partition gets finer:
\be 
var_p(X(0,t))=
\lim_{n \rightarrow \infty} var_p(X(0,t),\Pi_n)
\ee
We have not specified the nature of the convergence in the above 
definition, nor have we given any conditions on the sequence of partitions. 
In our simulations of the SAW and SLE we will use a sequence of 
uniform partitions for the $\Pi_n$. Changing the parameterization of 
the process is equivalent to changing the partitions used. As we will see, 
changing the parameterization can change the value of the $p$-variation. 

It is important to note that we do not take a supremum over all partitions
in the definition of the $p$-variation. If one takes a supremum one 
obtains a different quantity which, unfortunately, is also called the 
$p$-variation in the literature. (It is sometimes called the strong 
$p$-variation.)
For Brownian motion and $p=2$ the variation we have defined is the quadratic
variation studied by L\'evy \cite{levy}. If the sequence of partitions
satisfies $||\Pi_n|| = o(1/\log(n))$, then the quadratic variation of Brownian 
motion converges with probability one to $t$ \cite{dud}. But
if the condition is relaxed to  $||\Pi_n|| = O(1/\log(n))$ this 
is not true \cite{vega}. 
For other values of $p$, the $p$-variation has been shown to exist for
certain Gaussian processes including fractional Brownian motion and the 
local times associated with a symmetric stable L\'evy process. 
See the recent references \cite{fg, marrosa,  marrosb, shao} 
and references therein.

Now consider the SAW. If the scaling limit exists then 
$|\omega(t_j)-\omega(t_{j-1})|$ is of size $(t_j-t_{j-1})^\nu$. 
So it is natural to consider the $p$-variation with $p=1/\nu$.
We conjecture that 
\be 
var_{1/\nu}(\omega(0,t))=c \, t
\ee
for some constant $c$. 
If the $1/\nu$ variation exists for the scaling limit of the SAW, 
then it is trivial to see that 
\be 
E [var_{1/\nu}(\omega(0,t))]=c \, t
\ee
for some constant $c$. So the content of the conjecture is that the 
$1/\nu$ variation exists and is non-random. 
For the remainder of the paper, $p$ will be $1/\nu=4/3$, and so we will
drop the subscript $1/\nu$ on $var$.  

\begin{figure}[tbh]
\includegraphics{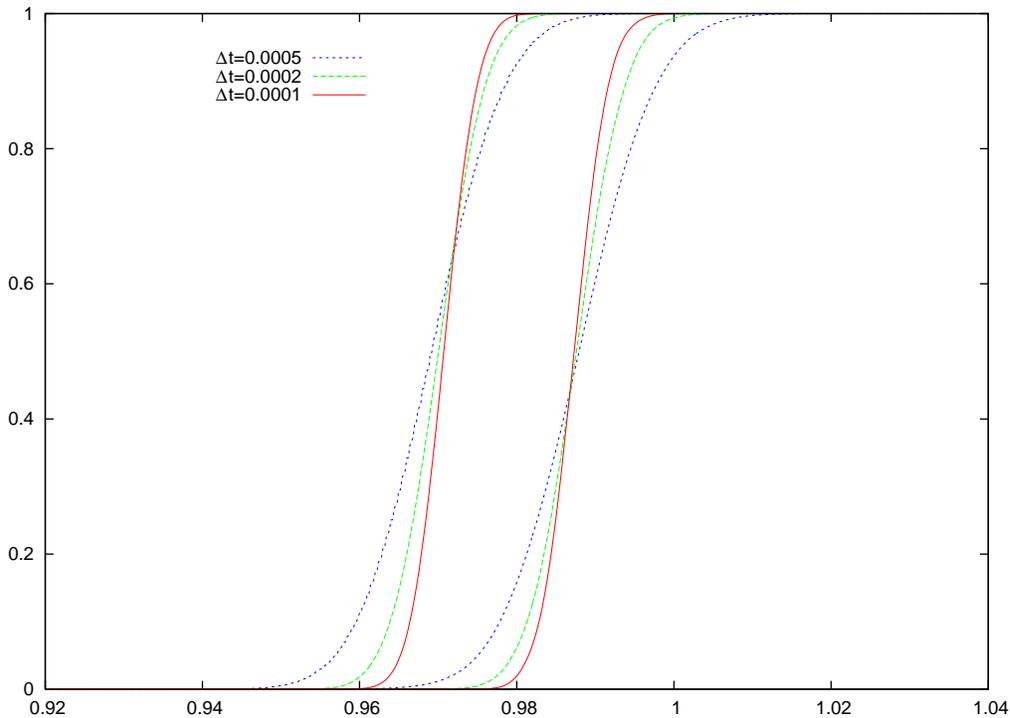}
\caption{ The set of curves on the right are the distributions of 
$var(\omega(0,1),\Delta t)$; those on the left  are the distributions of 
$var_{no}(\omega(0,1),\Delta t)$.
}
\label{label_saw_fvar}
\end{figure}

We support the conjecture with a simulation. We only consider 
a uniform partition with intervals of width $\Delta t$, and we take $t=1$. 
In this case we denote $var(\omega(0,1),\Pi)$ by $var(\omega(0,1),\Delta t)$. 
The distribution of this random variable 
is shown in figure \ref{label_saw_fvar} for several values of $\Delta t$. 
Two different sets of curves are shown in this figure. 
The curves on the right side of the figure are for the variation
defined above. 
Note that the scale on the horizontal axis begins at $0.92$, not $0$.
The figure clearly indicates that the distribution is converging to 
a step function as it should if the $1/\nu$ variation is constant.

\begin{figure}[tbh]
\includegraphics{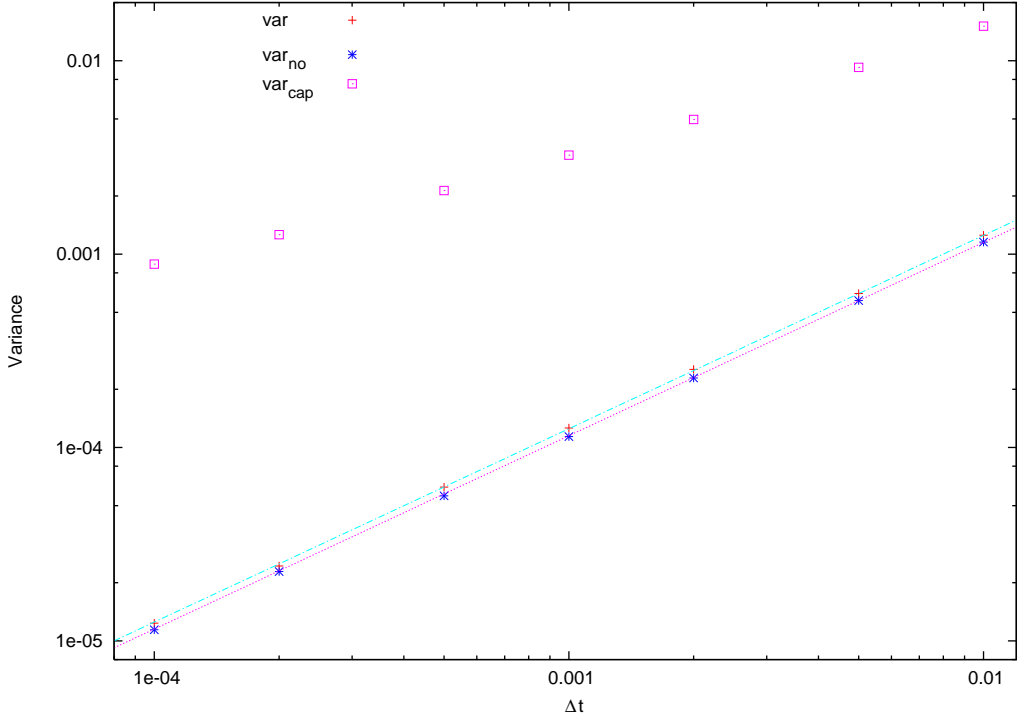}
\caption{ The variance of the three different definitions of the 
$p$-variation as functions of $\Delta t$. 
}
\label{label_saw_fvar_variance}
\end{figure}

To study the convergence quantitatively, in figure 
\ref{label_saw_fvar_variance} we plot the variance 
of the random variable $var(\omega(0,1),\Delta t)$ 
as a function of $\Delta t$.
Three curves are shown in this figure. 
The variance of $var(\omega(0,1),\Delta t)$ is the line of points 
labeled ``var.''
In this log-log plot the data is very well fit by a line with slope $1$. 
This corresponds to the variance being proportional to $\Delta t$. 
Note that this is what we would find if the process had independent, 
stationary increments. We do not expect the increments of the SAW to 
be independent, but this result suggests that the lengths of these increments
are weakly correlated. 

The fractional variation defined above involves the parameterization of the 
SAW. A different parameterization could give a different value to this 
variation. (We will see an example of this later.) Another definition of 
the variation that does not depend on the parameterization is the 
following. Let $\Delta t>0$. We define times $t_i$ as follows. 
$t_0=0$. Given $t_i$, we let $t_{i+1}$ be the first time after $t_i$
such that $|\omega(t_{i+1})-\omega(t_i)|=(\Delta t)^\nu$. The variation 
over the time interval $[0,t]$ is then defined to be $n \Delta t$,
where $n$ is the largest index with $t_n \le t$. We denote this variation 
by $var_{no}(\omega(0,t),\Delta t)$. The subscript $no$ stands for 
``no parameterization,'' to emphasize that this definition does not
depend on how the curve is parameterized.

The curves on the left of figure \ref{label_saw_fvar} show the variation
$var_{no}(\omega(0,1),\Delta t)$ for several values of $\Delta t$. 
The figure shows this variation is also converging to a constant. 
It also shows that the constant $var_{no}(\omega[0,t])$ is slightly 
less than the constant $var(\omega[0,t])$.
The variance of $var_{no}(\omega[0,t],\Delta t)$ is the line of points 
labeled ``var$_{\rm no}$''
in figure \ref{label_saw_fvar_variance}.
This data is also very well fit by a line with slope 1.  

\begin{figure}[tbh]
\includegraphics{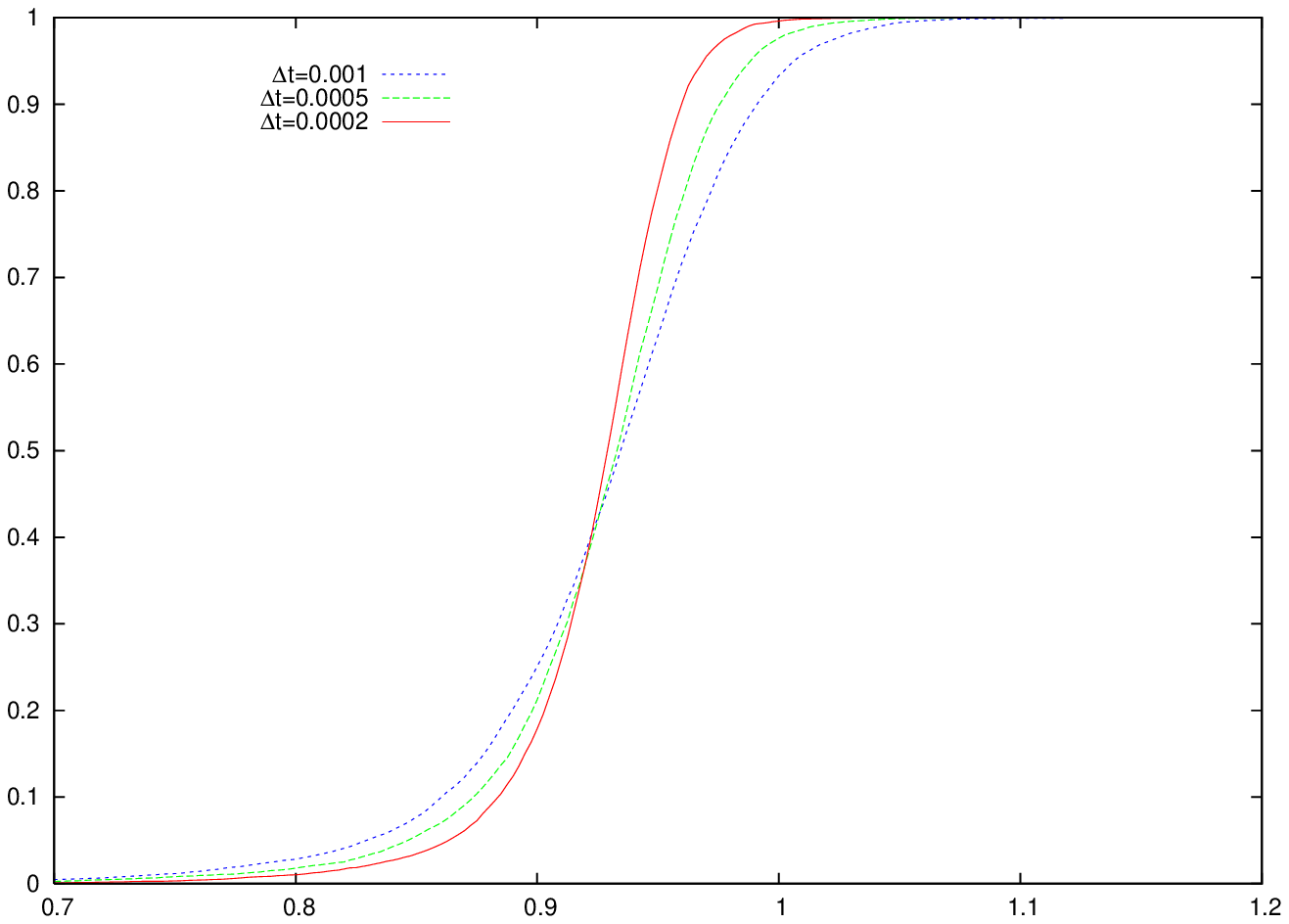}
\caption{The distribution of $var_{cap}(\omega(0,1),\Delta t)$, the random
variable that converges to the $p$-variation computed using the 
parameterization of the SAW by capacity. 
}
\label{label_saw_cp_fvar}
\end{figure}

We now return to our first definition of the $p$-variation, 
eq. \reff{fvar_def}, and consider what can happen if we use a different
parameterization. We are particularly interested in what happens when
we use the half-plane capacity to reparameterize the SAW. 
We denote the resulting variation by $var_{cap}(\omega(0,t),\Pi_n)$. 
We can numerically compute the capacity along a SAW and so we can simulate 
this random variable.
We have simulated it with a uniform partition. 
Computing the capacity is difficult and this limits the lengths of 
the SAW's we can study. This in turn restricts the simulations to 
larger values of $\Delta t$.
The distributions of this random variable for several choices of 
$\Delta t$ are shown in figure \ref{label_saw_cp_fvar}.
The first thing to note about this figure is the scale on the horizontal 
axis. These distributions have considerably larger variance than the 
curves in figure \ref{label_saw_fvar} with the same values of $\Delta t$. 
Nonetheless, it appears this this random variable is also converging 
to a constant as $\Delta t$ goes to zero. We also see that 
$var_{cap}(\omega(0,1))$ is slightly less than $var_{no}(\omega(0,1))$.
The variance of $var_{cap}(\omega(0,1),\Delta t)$ as a function of $\Delta t$ 
is the set of points labeled ``var$_{\rm cap}$''
in figure \ref{label_saw_fvar_variance}.
This data indicates the variance is converging to zero as $\Delta t$ goes
to zero, but it is not well fit by any line, at least for the range of 
$\Delta t$ shown.  

For each of the three variations we have studied, the plots of their
distributions for three values of $\Delta t$ very nearly have a common 
point of intersection. This provides a simple way to estimate the values
of these variations in the limit $\Delta t \rightarrow 0$. We
estimate 
\bea 
var(\omega(0,1)) &= 0.987 \pm 0.001 \nonumber \\
var_{no}(\omega(0,1)) &= 0.972 \pm 0.001 \nonumber \\
var_{cap}(\omega(0,1)) &= 0.92 \pm 0.01 
\eea
The error bars we have given are just guesses. As we discuss in 
section \ref{simulating} there are systematic errors in the simulations
which make it difficult to give reliable error bars. 

\bigskip

\section{Reparameterizing the SLE}
\label{reparamSLE}

If the scaling limit of the SAW is indeed SLE, then the results of the 
previous section suggest that each of the $1/\nu$ variations we considered 
should exist for SLE and give a parameterization which corresponds to 
the natural parameterization of the SAW. 
Of course, there is no way to compute the variation $var(\omega[0,t])$
for SLE since it requires knowing the parameterization we are trying to find. 
We must either use $var_{no}(\omega[0,t])$ or $var_{cap}(\omega[0,t])$. 
We use the former since its variance for a nonzero $\Delta t$ is smaller. 
Let $V_t$ be the random time on the SLE where 
\be
var_{no}(\gamma[0,V_t])=ct
\ee
where $c=var_{no}(\omega[0,1])$. 
We define $\hat{\gamma}(t)=\gamma(V_t)$ so that 
$var_{no}(\hat{\gamma}[0,t])=var_{no}(\omega[0,t])$
Then $\hat{\gamma}(t)$ and $\omega(t)$ 
should have the same distribution as parameterized curves. 
We test this by comparing the distributions of $\omega(1)$ and 
$\hat{\gamma}(1)=\gamma(V_1)$.

We must compute the constant $c=var_{no}(\omega[0,1])$ by simulation.
This requires computing the variation $var_{no}(\omega(0,1),\Delta t)$
for several $\Delta t$ and then estimating the limit as 
$\Delta t \rightarrow 0$. 
However, when we compute the variation $var_{no}(\gamma[0,V_t])$ we must 
also use a nonzero $\Delta t$ and then attempt to take a limit as 
$\Delta t \rightarrow 0$. These extrapolations have some error in them. 
We attempt to minimize this error by the following trick. 
For a given $\Delta t$ we define $V_t$ by 
\be
var_{no}(\gamma[0,V_t],\Delta t)= var_{no}(\omega[0,1],\Delta t) \, t
\ee
and then let $\hat{\gamma}(t)=\gamma(V_t)$. In other words, for the 
constant $c$ we use the estimate of $c$ that comes from the SAW simulation
using the same value of $\Delta t$ that we use to compute the $p$-variation
for the SLE. (Note that $V_t$ now has some small dependence on $\Delta t$.)

\bigskip

\begin{figure}[tbh]
\includegraphics{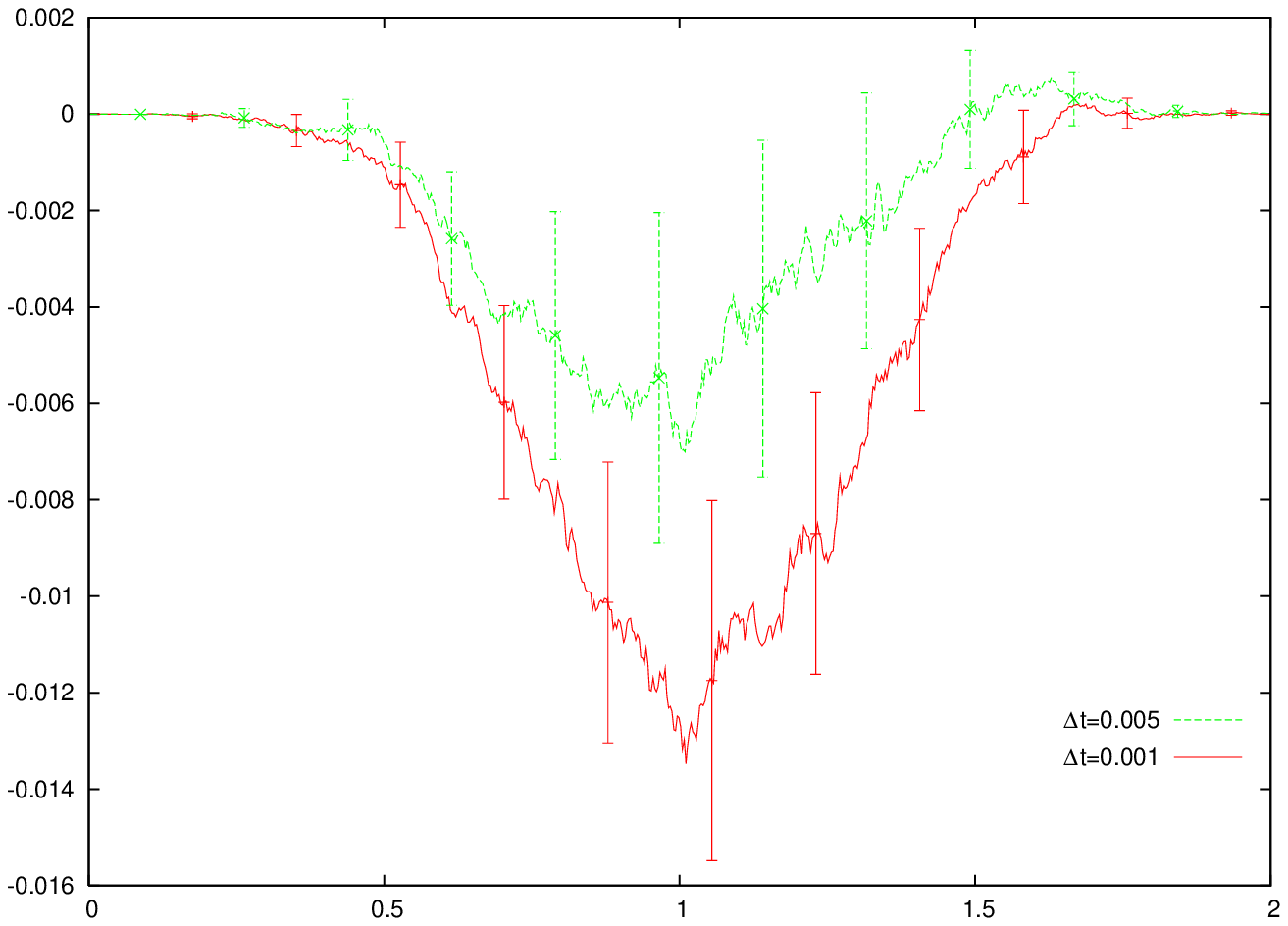} 
\caption{The difference of the distributions of $R$ for $\omega(1)$, 
the SAW with its natural parameterization, and $\hat{\gamma}(1)$, the 
SLE parameterized by the variation $var_{no}$. 
}
\label{saw_sle_fix_fvar_dist}
\end{figure}

\begin{figure}[tbh]
\includegraphics{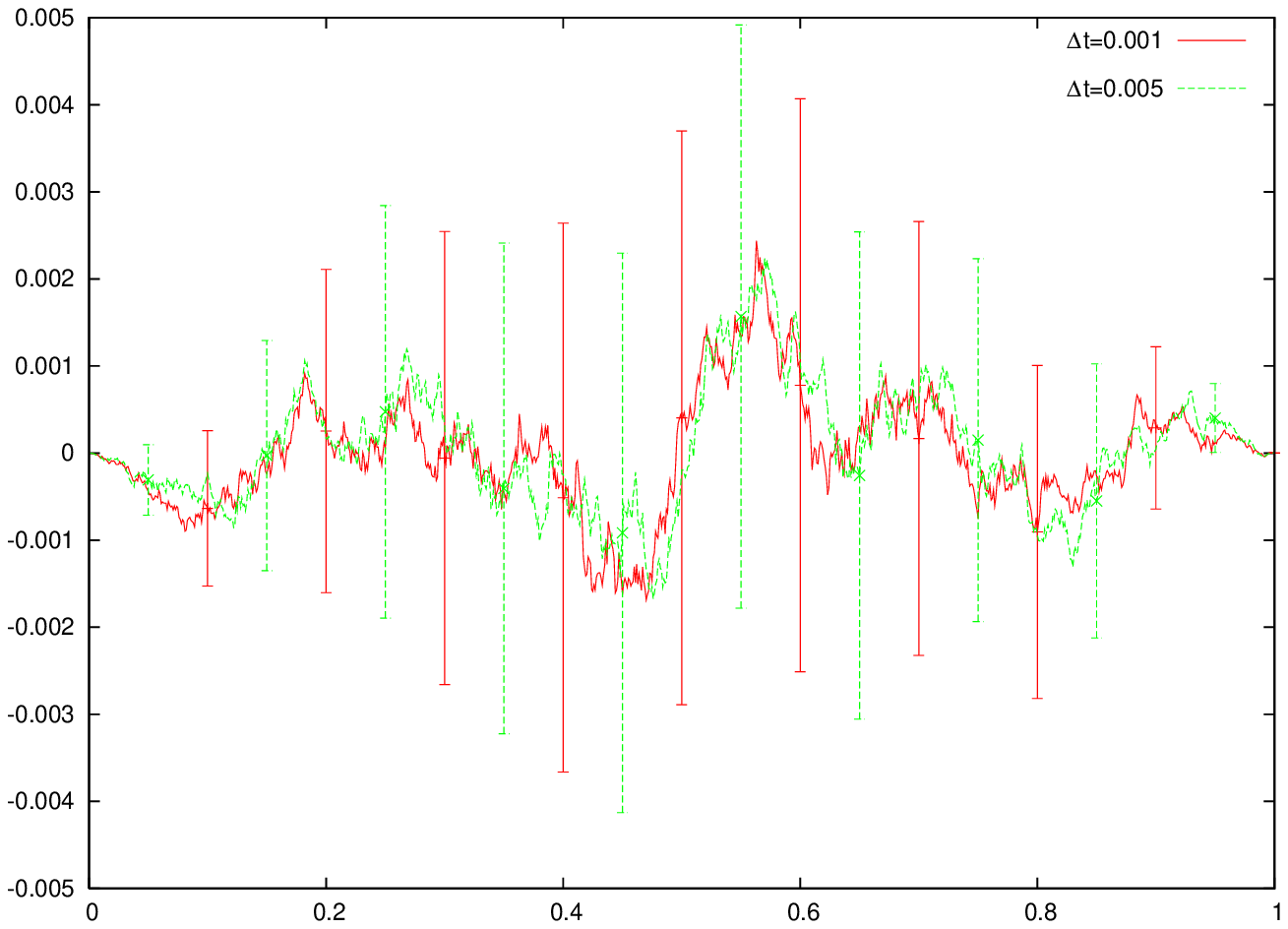}
\caption{The difference of the distributions of $\Theta$ for $\omega(1)$, 
the SAW with its natural parameterization, and $\hat{\gamma}(1)$, the 
SLE parameterized by the variation $var_{no}$. 
}
\label{saw_sle_fix_fvar_angle}
\end{figure}

As before we compute the distributions of the random variables $R$ 
and $\Theta$.
In figure  \ref{saw_sle_fix_fvar_dist} we show the difference between 
the distributions of $R$ for $\omega(1)$ and $\hat{\gamma}(1)$.
Two plots are shown corresponding to two different values of $\Delta t$. 
The difference is larger than that seen in figure \ref{saw_sle_fix_cap_dist}. 
However, the figure also shows that the difference depends strongly 
on $\Delta t$. 
In figure  \ref{saw_sle_fix_fvar_angle} we show the difference between 
the distributions of $\Theta$ for $\omega(1)$ and $\hat{\gamma}(1)$.
For the angle we see almost no dependence on the choice of $\Delta t$. 

\begin{figure}[tbh]
\includegraphics{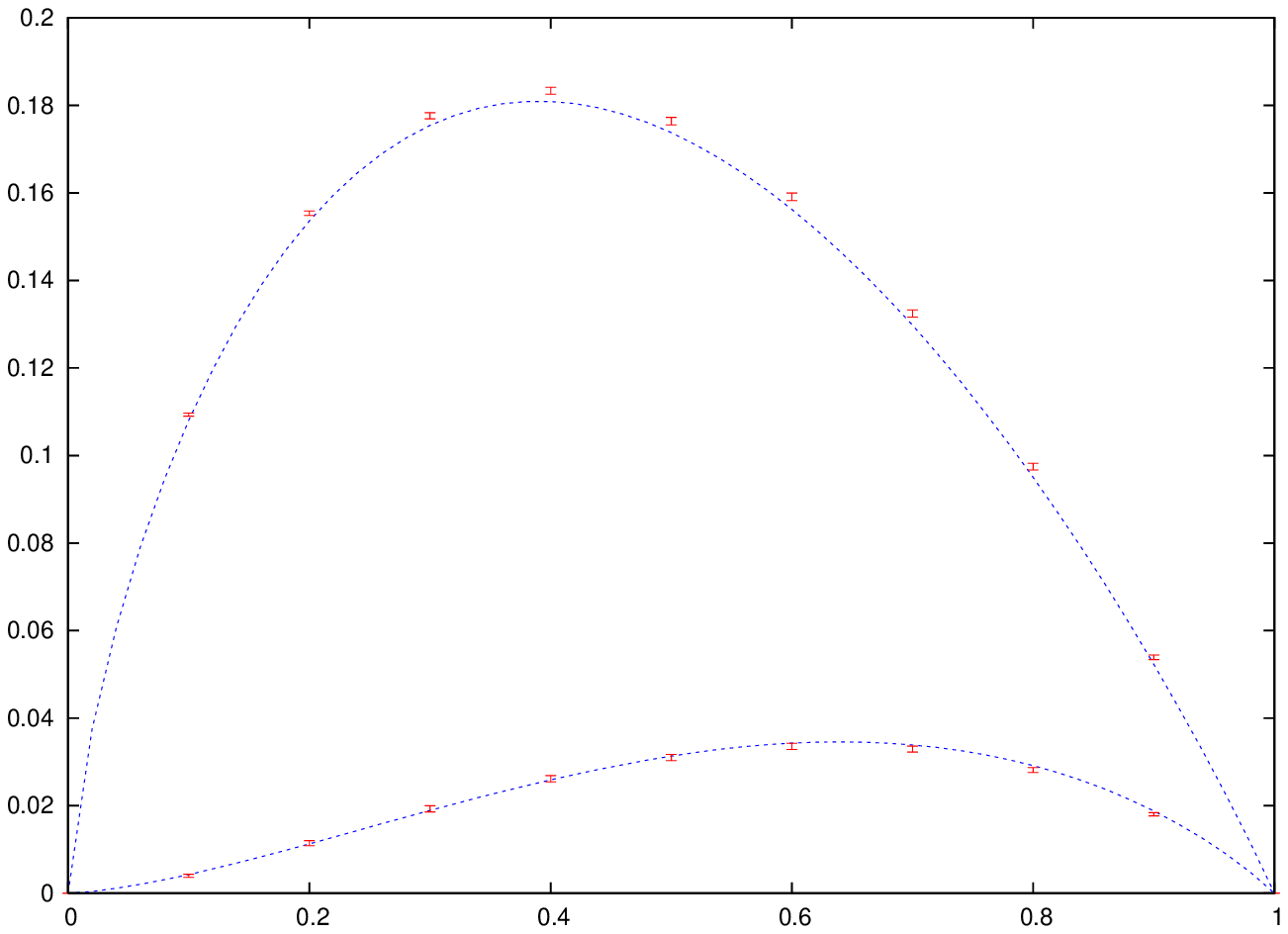}
\caption{ The covariance $\sigma_{ij}(t)$ for the SAW and the covariance
$\hat{\rho}_{ij}(t)$ for SLE using the variation $var_{no}$ 
for its parameterization.
The top curve is $i=j=2$, the bottom curve is $i=j=1$. 
The SLE covariance is plotted with a curve and no error bars, while 
individual points with error bars are plotted for the SAW. 
}
\label{saw_sle_fix_fvar_covar}
\end{figure}

Next we test the equivalence of $\hat{\gamma}(t)$ and $\omega(t)$ 
by comparing their covariances. As we saw before, the non-trivial 
part of these covariances is given by
\be
\hat{\rho}_{ij}(t)= E  \, [(\hat{\gamma}_i(1) 
- \hat{\gamma}_i(t)) \hat{\gamma}_j(t)] 
\ee
and
\be
\sigma_{ij}(t)= E  \, [(\omega_i(1) - \omega_i(t)) \omega_j(t)] 
\ee
These functions are shown in figure \ref{saw_sle_fix_fvar_covar} 
for $i=j=1$ and $i=j=2$. The SAW and SLE covariances agree reasonably well. 
As in figure \ref{saw_sle_fix_fvar_dist} there is significant 
dependence on the $\Delta t$ used to compute the $p$-variation of the SLE.
For some values of $\Delta t$ the agreement between the two covariances 
is not as good as that shown in the figure. 

It is tempting to compare the plots in figures 
\ref{saw_sle_fix_cap_covar} and \ref{saw_sle_fix_fvar_covar}.
The first thing to keep in mind is that in figure 
\ref{saw_sle_fix_cap_covar}, $t=1$ corresponds to $hcap=2$ while 
in figure \ref{saw_sle_fix_fvar_covar}, $t=1$ corresponds to 
$var_{no} = 1$. Simulations show that the average capacity for a SAW
with $var_{no} = 1$ is close to $0.5$. Changing the capacity by a factor
of $4$ is the same as changing spatial scales by a factor of $2$, and 
so corresponds to a change in the covariance by a factor of $4$. 
Indeed the vertical scales in the two figures differ by a factor of
approximately $4$. 
However, we should emphasize that 
the two plots use completely different parameterizations
of the processes and so comparing them is not particularly meaningful. 
These two plots should be different even if they are rescaled to have the 
same vertical scale. 

\section{Simulating SAW and SLE}
\label{simulating}

One method for simulating SLE is to approximate the conformal
map $g_t$ by the composition of a large number of conformal maps
which are known explicitly. 
A particular implementation was studied by  R. Bauer \cite{bauer}.
This method is also briefly discussed for the radial case in 
\cite{mr}. Here we give a brief overview. 
A more detailed review may be found in  \cite{tk_sle}. 

We consider the time interval $[0,1]$ and 
divide it into $N$ subintervals.
(We could use a uniform partition of $[0,1]$, but one obtains a more 
uniform distribution of points along the SLE curve by using a non-uniform 
partition of the time interval. See \cite{tk_sle} for details.)
The conformal map $g_t$ can be written as the composition of 
$N$ conformal maps, each of which corresponds to the solution of 
the SLE equation \reff{sle_eq} over one of the small time intervals. 
Each of these $N$ conformal maps is approximated by a simple 
conformal map. We use the conformal map that takes the half plane 
minus a slit onto the half plane. The parameters for these conformal 
maps (e.g., the length and angle of the slit) are chosen so that 
the conformal maps have the correct capacity and so that the 
driving function corresponding to the composition of these conformal maps 
approximates the original driving function $\sqrt{\kappa} B_t$. 
For example, this can be done so that the driving function of the 
approximation agrees with $\sqrt{\kappa} B_t$ at the times which are 
endpoints of the $N$ subintervals. 

To compute a point on the SLE trace, one must apply $O(N)$ of these
conformal maps. So the time to compute a single point is $O(N)$.
If one wants to compute $N$ points on the SLE trace this will 
take a time $O(N^2)$. The paper \cite{tk_sle} shows how to implement 
this method of simulating the SLE so that the time required to 
compute a single point is approximately $O(N^{0.4})$ rather than $O(N)$. 

Certain features of the SLE are relatively easy to study by 
simulation. For example, to compute the distribution of the SLE trace
at a fixed time, one need only compute a single point on the SLE trace.
One can use a relatively large value of $N$ and generate 
a large number of samples. 
By contrast, computing the $1/\nu$ variation of the curve 
involves a double limit. One needs to let $N \rightarrow \infty$ and 
then let the time interval, $\Delta t$, used to compute the 
variation go to zero. In practice this means that one must take a 
small $\Delta t$ and then take $N$ large enough that $1/N$ is small 
compared to $\Delta t$. Thus the simulations that require computing 
the variation are among the most difficult. 

We simulate the SAW using the pivot algorithm
\cite{ms}. The particular implementation of the pivot algorithm we use 
is found in \cite{tk_pivot}. This algorithm is fast in two dimensions, and 
one can study quantities like the location of the SAW at a fixed time
for walks with a million steps. The difficult part of the SAW 
simulations is computing the capacity of a walk. This can only be 
done for much shorter walks.

We compute the capacity of a SAW using the zipper algorithm  \cite{kuh,mr}. 
We give a brief explanation of how the algorithm
works in our context and refer the reader to \cite{mr} for more 
detail. Consider a curve $\omega(t)$ in the upper half plane which 
starts at the origin and ends at $P$, e.g., a self-avoiding walk. 
To compute the capacity we need to find the conformal map $g(z)$ that 
takes the half plane $\half$ minus the curve onto $\half$. It should 
be normalized so that $g(\infty)=\infty$ and $g^\prime(\infty)=\infty$. 
This determines the map up to a translation. 
The capacity is the coefficient of $1/z$ in the expansion of $g(z)$ 
about $\infty$. It does not depend on the choice of translation. 

Let $0=z_0,z_1,z_2, \cdots, z_n=P$ be points along the curve.
Let $C_i$ denote the portion of the curve from $z_{i-1}$ to $z_i$. 
Let $g_1(z)$ be the conformal map that takes $\half \setminus C_1$ 
onto $\half$. In addition to the above normalizations, 
we require $g_1(z_1)=0$. Then $g_1$ maps $C_2 \cup \cdots \cup C_n$ 
to a curve that starts at $0=g_1(z_1)$ and ends at $g_1(z_n)$ and 
passes through the points $g_1(z_i)$, $i=2,\cdots, n-1$. 
Now we let $g_2(z)$ be the conformal map that takes 
$\half \setminus g_1(C_2)$ onto $\half$. Then $g_2 \compose g_1(z)$ takes 
$C_3 \cup \cdots \cup C_n$ to a curve that starts at 
$0=g_2 \compose g_1(z_2)$, ends at $g_2 \compose g_1(z_n)$ 
and passes through the points $g_2 \compose g_1(z_i)$, $i=3,\cdots, n-1$. 
We continue to remove one $C_i$ at a
time. $g_i(z)$ is the conformal map that takes 
$\half \setminus g_{i-1} \compose \cdots \compose g_2 \compose g_1(C_i)$ 
onto $\half$. 
The composition $g(z)=g_n \compose \cdots \compose g_2 \compose g_1(z)$ 
is then the desired conformal map. 
The capacity of $g$ is the sum of the capacities of the $g_i$.

Note that even if the original curve is piece-wise linear, the 
segments of the form $g_1(C_1)$, $g_2 \compose g_1(C_2) \cdots$ are not. 
The key idea of the zipper algorithm is to approximate $g_i$ by 
the conformal map $h_i$ 
that takes $\half \setminus D_i$ onto $\half$ where $D_i$ is 
a curve that starts at $0$ and ends at 
$g_{i-1} \compose \cdots g_2 \compose g_1(z_i)$.
So $D_i$ has the same endpoints as 
$g_{i-1} \compose \cdots \compose g_2 \compose g_1(C_i)$.
The curve $D_i$ is chosen so that the map $h_i$ is relatively easy to compute. 
One possibility is to take $D_i$ to be a line segment from $0$ to 
$g_{i-1} \compose \cdots \compose g_2 \compose g_1(z_i)$. 
Then there is an explicit expression 
for the inverse of $h_i$ and $h_i$ itself may be found by a Newton's
method. (Some care is needed in implementing Newton's method \cite{mr}.)
Another choice is to take $D_i$ to be the arc of the circle that is 
perpendicular to the real axis and connects $0$ and 
$g_{i-1} \compose \cdots \compose g_2 \compose g_1(z_i)$. 
This has the advantage that $h_i$ 
only involves linear fractional transformations and a square root. 

Because of the need to compute 
$g_{i-1} \compose \cdots \compose g_2 \compose g_1(z_i)$ for every point
$z_i$, the zipper algorithm requires a time $O(n^2)$. The idea in 
\cite{tk_sle} that speeds up the SLE simulation may be applied here
to dramatically reduce the time required by this algorithm.

The curve $D_i$ starts and ends at the same points as 
$g_{i-1} \compose \cdots \compose g_2 \compose g_1(C_i)$, 
but does not necessarily approximate it well. For a fixed curve one 
improves the approximation to the capacity by using more points
$z_1, \cdots, z_n$ along the curve. We do not do this. Thus the capacity
we compute is not exactly the capacity of the SAW. However, it is 
the exact capacity (up to numerical errors) 
of some curve that passes through the same lattice points as the SAW.
We expect that this curve will have the same scaling limit as the SAW.

As in any Monte Carlo simulation there is error because we do not run
the simulation forever. We refer to this as ``statistical error.''
It is relatively easy to estimate (we use batched means). 
There are also systematic errors in the simulations which we now discuss.

For the SAW, the scaling limit is given by a double limit. First we should 
let the length of the walk go to infinity and then take the lattice 
spacing to zero. In practice we fix a large $N$ and simulate SAW's 
with $N$ steps but only use the first $N^\prime$ steps where
$N^\prime$ is significantly smaller than $N$. 
The simulation is done with a unit lattice, but then we rescale the 
walks by a factor of $(N^\prime)^{-\nu}$. So the portion of the 
rescaled SAW that we use has a 
size of order $1$ and approximates the scaling limit of the SAW 
with the natural parameterization running from $t=0$ to $t=1$. 
For the SAW simulations we will give the 
number of iterations of the Monte Carlo Markov Chain. This number is 
typically very large, but the samples generated are highly correlated.

For the SLE, $N$ denotes the number of conformal maps used in the 
approximation. 
For the SLE simulations we will give the number of SLE's generated.
These samples are independent. 

We end this section with a discussion of the parameters used in the 
various simulations.
In figures \ref{saw_sle_fix_time_dist} and \ref{saw_sle_fix_time_angle},
the SAW is simulated using $N=1,000,000$ and $N^\prime=200,000$. 
The simulation ran for 1 billion iterations. 
The SLE simulation used $N=100,000$, and 1 million samples were generated. 
These large numbers are possible because all we need to compute for the 
SLE is the location of its position at $t=1$.

In figures \ref{saw_sle_fix_cap_dist} and \ref{saw_sle_fix_cap_angle}
the SLE simulation used is the same as that in figures 
\ref{saw_sle_fix_time_dist} and \ref{saw_sle_fix_time_angle}.
The SAW simulation requires computing the capacity of the walk. 
This limits the simulation to significantly shorter walks. 
We take $N=100,000$ and rescale the walk by a factor of $N^{-\nu}$. 
For this rescaled walk we compute
the capacity along the walk up until $hcap=0.1$. (The number of 
steps at which this occurs is random but well below $N$. The mean number 
of steps at which the capacity is $0.1$ is about a third of the total number
of steps.) We then rescale the SAW again by a factor of $\sqrt{20}$ 
so that the random time we are finding is where the capacity is $2$. 
As we have noted the samples generated by the SAW simulation are highly 
correlated. If we are studying an observable which is trivial to compute,
we might as well compute the value of the observable at every time step.
The samples will be highly correlated, but there is very little cost in 
terms of computation time. 
For observables that are not trivial to compute, e.g., observables that 
involve computing the capacity, if we computed the observable at every 
time step we would spend most of the computation time on computing the 
observable. For such observables we only compute it 
every $s$ time steps. In this SAW simulation 
we took $s=100,000$. The simulation was run for 8 billion iterations, so 
$80,000$ samples of the observables were computed. 

For the first covariance plot, figure \ref{saw_sle_fix_cap_covar},
the SLE simulation is done with $N=10,000$. We compute the SLE covariance 
at $50$ equally spaced times. 200,000 samples of the SLE were 
generated. For the SAW simulation, we take $N=100,000$, rescale the 
walk by a factor of $N^{-\nu}$, and then compute
the capacity along the walk up until $hcap=0.1$. 
We then rescale the SAW again by a factor of $\sqrt{20}$ 
so that the random time we are finding is where the capacity is $2$. 
Computing the capacity of the SAW is slow, so we only compute it 
every $100,000$ iterations of the Markov chain. We ran the chain for 
$8$ billion iterations, so we generated $80,000$ samples of the SAW 
covariance. 

There are two simulations of the SAW in figure \ref{label_saw_fvar}.
In both of them $N=1,000,000$ and we use $N^\prime=500,000$ and
$N^\prime=1,000,000$. For most observables, taking $N^\prime$ at or near
$N$ is a bad idea - it will change the distribution of the observable. 
However, for the two variations shown in the figure we find no difference 
between the variations computed using $N^\prime=500,000$ 
and $N^\prime=1,000,000$. The curves shown use the latter value. 
Both simulations were run for $1$ billion iterations. 

In the simulation of the SAW used for figure \ref{label_saw_cp_fvar}, 
we take $N=100,000$ and $N^\prime=10,000$. The relatively small value
of $N^\prime$ is because of the difficulty in computing the capacity of a SAW. 
We run the Markov chain for $1$ billion iterations, but only 
compute $var_{cap}$ every $100,000$ iterations for a total of 
only $10,000$ samples.

In figures \ref{saw_sle_fix_fvar_dist} and \ref{saw_sle_fix_fvar_angle}
the SAW simulation used is the same as that in figures 
\ref{saw_sle_fix_time_dist} and \ref{saw_sle_fix_time_angle}.
For the SLE simulation we took $N=500,000$. However, we only 
computed every fifth point along the SLE, so we compute a total of $100,000$
points on the SLE. We generated 146,000 samples. This is the longest
simulation in this paper. It required approximately 150 cpu-days.
The time needed for this simulation is significantly reduced by the 
following trick. We compute the variation $var_{no}$ up to the time
when it is $1$. This time on the SLE when this is attained is random,
but occurs well before the end of the SLE we are computing. So we 
actually only need to compute some initial fraction of the $100,000$ points on 
the SLE. 

Finally, for figure \ref{saw_sle_fix_fvar_covar} the SAW simulation 
uses $N=1,000,000$ and $N^\prime=400,000$. The SAW simulation was 
run for $1$ billion iterations. For the SLE simulation, $N=250,000$
but only every fifth point was computed. Again, as discussed in the 
preceding paragraph we need only compute an initial portion of the SLE.
127,000 samples of the SLE were generated. 

\section{Conclusions}
\label{conclusion}

The main conclusion of this paper is that the $p$-variation with 
$p=1/\nu$ of the SLE trace provides a parameterization that corresponds to 
the natural parameterization of the SAW. Simulations show that the 
$p$-variation of the SAW is non-random. A trivial scaling argument then 
implies it is proportional to the natural parameterization of the SAW.
Thus the SLE with the $p$-variation as its parameterization and the 
SAW with its natural parameterization should agree as parameterized curves. 
Two tests were done to check this agreement. 
The distributions of the random curves at a fixed time were compared, and 
the covariances of the two processes were compared. Good agreement 
was found. 
 
A secondary conclusion of this paper concerns using the half-plane capacity
to reparameterize the SAW. With this reparameterization it should agree with 
the SLE as parameterized curves.
The same two tests of their equivalence as parameterized curves 
showed good agreement.

For random fractal curves like the SAW or the loop-erased random walk
the exponent $\nu$ should be equal to the reciprocal of the Hausdorff 
dimension of the curve. Beffara \cite{befa,befb} has proved that 
if $\gamma$ is the SLE trace with parameter $\kappa$, 
then the Hausdorff dimension of $\gamma$ is $1+\kappa/8$ a.s.
This suggests that the $p$-variation of SLE should exist for 
$p=\nu^{-1}=1+\kappa/8$.
A natural question is to consider this variation for other discrete models
and see if it exists and is non-random as it appears to be for the SAW.
Preliminary simulations of the loop-erased random walk indicate that its
$1/\nu$ variation exists and is non random for $\nu=4/5$. 
This was one of the models considered by Schramm in his original paper
\cite{schramm}. Lawler, Schramm and Werner \cite{lsw_lerw} 
have proved that its scaling limit converges to SLE with $\kappa=2$.

\bigskip
\bigskip

\noindent {\bf Acknowledgments:}
The Banff International Research Station made possible many useful
interactions. In particular, the author thanks 
David Brydges, Greg Lawler, Don Marshall, Daniel Meyer, Yuval Peres, 
Stephen Rohde, Oded Schramm, Wendelin Werner and Peter Young for 
useful discussions. 
This work was supported by the National Science Foundation (DMS-0201566
and DMS-0501168).

\end{document}